\documentclass[preprint,12pt]{elsarticle}
\usepackage{amssymb}
\usepackage{comment}
\usepackage{multirow}
\usepackage{multicol}
\usepackage{adjustbox}
\usepackage{xcolor}
\usepackage{booktabs}
\usepackage{multirow}
\usepackage{comment}
\usepackage[colorlinks]{hyperref}
\usepackage{adjustbox}
\usepackage{amsmath}
\usepackage{amssymb}
\usepackage{amsthm}
\usepackage{ dsfont }
\usepackage{graphicx}
\usepackage{float}
\usepackage{pdfpages}
\usepackage{etoolbox}
\usepackage{lipsum}
\usepackage[short]{optidef}
\usepackage{xcolor}
\usepackage{caption}
\usepackage{fancyhdr,graphicx,amsmath,amssymb}
\usepackage[linesnumbered,ruled,vlined,noend]{algorithm2e}  
\usepackage{float}
\usepackage{geometry}
\newgeometry{left=2cm,right=2cm,top=2cm,bottom=2cm}
\journal{Logistics, MDPI}
\begin{document}
\begin{frontmatter}

\title{A Mixed-integer Linear Formulation for Dynamic Modified Stochastic p-Median Problem in a Competitive Supply Chain Network Design}

\author[inst1]{Amir Hossein Sadeghi\corref{mycorrespondingauthor}}
\cortext[mycorrespondingauthor]{Corresponding author}
\ead{asadegh3@ncsu.edu}
\author[inst2]{Ziyuan Sun}
\author[inst1]{Amirreza Sahebi Fakhrabad}
\author[inst3]{Hamid Arzani}
\author[inst4]{Robert B. Handfield}

\affiliation[inst1]{
            organization={Department of Industrial and Systems Engineering}, 
            addressline={North Carolina State University},
            state={NC},
            country={USA}}

\affiliation[inst2]{
            organization={Department of Mechanical, Industrial and Aerospace Engineering},
            addressline={Concordia University},
            state={Québec},
            country={Canada}}
       
\affiliation[inst3]{
            organization={Rotman School of Management},
            addressline={University of Toronto},
            state={Ontario},
            country={Canada}}
            
\affiliation[inst4]{ 
            organization = {Department of Business Management Poole College of Management},
            addressline={North Carolina State University},            
            state={NC},
            country={USA}}

\begin{abstract}
The Dynamic Modified Stochastic p-Median Problem (DMS-p-MP) is an important problem in supply chain network design, as it deals with the optimal location of facilities and the allocation of demand in a dynamic and uncertain environment. In this research paper, we propose a mixed-integer linear formulation for the DMS-p-MP, which captures the key features of the problem and allows for efficient solution methods. The DMS-p-MP adds two key features to the classical problem: (1) it considers the dynamic nature of the problem, where the demand is uncertain and changes over time, and (2) it allows for the modification of the facility locations over time, subject to a fixed number of modifications. The proposed model is using robust optimization in order to address the uncertainty of demand by allowing for the optimization of solutions that are not overly sensitive to small changes in the data or parameters. To manage the computational challenges presented by large-scale DMS-p-MP networks, a Lagrangian relaxation (LR) algorithm is employed. Our computational study in a real-life case study demonstrates the effectiveness of the proposed formulation in solving the DMS p-Median Problem. According to the results, the number of opened and closed buildings remains unchanged as the time horizon increases. This is due to the periodic nature of our demand. This formulation can be applied to real-world problems, providing decision-makers with an effective tool to optimize their supply chain network design in a dynamic and uncertain environment.
\end{abstract}

\begin{keyword}
p-Median Problem; Supply Chain Network Design; Dynamic Allocation; Robust Optimization; Lagrangian Relaxation
\end{keyword}

\end{frontmatter}

\section{Introduction}
Mobile grocery stores, also known as grocery trucks (GT) or pop-up grocery stores, have emerged as a new trend in the retail industry. GTs typically stock a wide variety of fresh fruits and vegetables, as well as other grocery items. They are usually equipped with refrigeration units to keep the food fresh and have a payment system in place for customers to buy the products.

These mobile stores bring convenience and accessibility to customers in under-served or rural areas, or in areas where traditional brick-and-mortar grocery stores are not present \cite{fahlevi2019leadership,nader2018contrarian}. The trend of grocery trucks is expected to continue to grow in the coming years as more and more people look for ways to improve access to healthy food in their communities. In the U.S., the mobile food vendors market is valued at 1.16 billion dollars in 2021 and is predicted to grow at a rate of 6.4\% annually from 2022 to 2030, driven by the increasing trend of culinary arts and the preference of young people for different dining experiences over the traditional dining in restaurants \cite{noauthor_us_nodate}.

GTs can also be used to provide food to people in emergency situations, like natural disasters or power outages \cite{hecht2019urban,pelling2001natural}, as well as address food insecurity by increasing access to healthy food options. GTs are valuable resources for alleviating food insecurity, and these trucks are often operated by non-profit organizations, local governments, and community groups that want to address food insecurity and promote healthy eating \cite{mohan2013improving}. They can help to address food insecurity in several ways: 1. Accessibility: Grocery trucks can bring fresh produce and other grocery items to low-income or rural areas where there is limited access to traditional grocery stores and supermarkets. This makes it easier for residents in these areas to access healthy food options; 2. Affordability: Many grocery trucks accept SNAP (Supplemental Nutrition Assistance Program) benefits and other forms of government assistance, making it more affordable for low-income families to purchase healthy food; 3. Education: Some grocery trucks provide educational programs and cooking demonstrations that teach customers how to prepare healthy meals. This can help to improve the overall health and nutrition of the community; 4. Variety: Grocery trucks can also offer a wider variety of fresh produce and other grocery items than traditional corner stores or convenience stores, which may only carry a limited selection of products.

Overall, grocery trucks are a creative and innovative way to address food insecurity and improve access to healthy food in communities that need it the most. However, the success of a mobile grocery store is heavily dependent on the location where it is parked \cite{esparza2014trade, wessel2012place}. Various factors should be considered when selecting a location for a grocery/food truck, such as dynamic demand, visibility, accessibility, and competition. 

Logistics for mobile grocery stores refer to the process of planning, coordinating, and controlling the movement and storage of goods, services and information from the point of origin to the point of consumption. This includes transportation, inventory management, warehousing, and distribution of products. For mobile grocery stores, logistics also includes the planning and coordination of the routes and schedule of the mobile store, as well as the management of the supply chain and inventory. This includes sourcing products from suppliers, managing inventory levels and restocking the mobile store as needed, and coordinating delivery schedules.

Effective logistics management is crucial for the success of mobile grocery stores, as it can impact delivery times, cost, and overall customer satisfaction. By optimizing logistics, mobile grocery store owners and operators can improve their business by reducing costs and increasing efficiency, ultimately resulting in better customer service \cite{restuputri2022customer,bourlakis2008food}. In this regard, we first provide a comprehensive literature review of the proposed location-allocation models in the literature in section \ref{sec:lit}, then provide a mathematical formulation for our model considering the uncertainty of demand in section \ref{sec:model}. Sections \ref{sec:example} and \ref{sec:conclusion} show the application of our proposed model in the real-world mobile grocery location problem and provide insights, respectively.

\section{Literature Review}\label{sec:lit}
Location science acts on a significant role in modern development in different disciplines, including business, economics, computer science, geography, military, transportation, etc.  With historical records, the first optimal location problem was proposed by Pierre de Fermat to find the geometric median among three points. The problem was solved by Evangelista Torricelli soon after. This is undoubtedly logical that the single facility minimum Euclidean distance problem is named Fermat-Torricelli problem (Known as the most famous Weber problem, the distances are looked upon as the weights of nodes). \cite{laporte2019introduction} is a book that comprehensively introduces the development and forecast of the subject of location science.

The beginning of the new era of location science was around the 1960s. \cite{berge1957two} investigates the problem of finding the minimum coverage on a graph originally. \cite{miehle1958link, cooper1963location} raises the best known $p$-median problem, which is to pick the exact $p$ of facilities to open that minimizes the total transportation cost. Followed by, \cite{hakimi1964optimum, hakimi1965optimum} releases Hakimi’s node optimally theorem to demonstrate that the optimal solution on a continuous graph for absolute median problems is always on the graph’s vertices, which means many network problems can simply transfer to discrete version problems. \cite{hakimi1964optimum} introduces the concept of absolute center to find the decision of satisfying the mini-max. As for solving more complex realistic optimization tasks in the next decade, mixed-integer linear programming (MILP) was widely used to deal with location problems. \cite{balinski1965integer, revelle1970central, toregas1971location} are early works to use MILP to solve the uncapacitated facility location problem, discrete $p$-median problem, and covering-location problem, respectively.

The deterministic model means all model parameters are known with certainty; nevertheless, unpredictability or randomness always occurs in real life. Stochastic programming or robust optimization is the solution to reduce the impact of uncertainties in predicting. Researchers have been quite interested in topics to combine location models with stochastic programming methods in the last few decades. \cite{louveaux1986discrete, louveaux1992dual} modify the uncapacitated facility location and p-median problems with stochastic variables, including the customer's demand, the unit cost of satisfying the customer from a specific location, and the fixed cost of locating at the candidate location. \cite{laporte1994exact} assumes the customer's demand is stochastic for the capacitated facility location problem. \cite{current1998dynamic} studies the number of opening facilities is uncertain. For improving customer service, based on $p$ initial opening facilities, \cite{berman2008p} allows increasing $r$ new opening facilities depending on the customer's demand, where $0\leq r\leq q$. \cite{sonmez2012decomposition} is the opposite of the previous article, $r$ existing facilities have to be closed, where $0\leq r\leq q \leq p$. \cite{liu2022employing} considers the covering-type location problems under demand uncertainties. \cite{snyder2006facility, correia2019facility} are overviews of facility location under uncertainty.

In reality, the vast majority of practical location problems can take into consideration time. The setting of the problem does not rely on time, called `single-period' or `static'; conversely, it is called `multi-period' or `dynamic' if the decision varies with different parameters at each period. \cite{ballou1968dynamic, sweeney1976improved} are early works about locating and re-locating a single facility within a time span. \cite{scott1971dynamic} continuous the previous works for locating multiple facilities. \cite{wesolowsky1973dynamic} merges the conception of the Weber problem with the multi-period facility location problem to find the best opening facility in each time period. \cite{warszawski1973multi, cavalier1985sequential, drezner1995dynamic, hakimi1999locations} studies the dynamic uncapacitated facility location problem, location problems on networks, network p-median problems, and network center problems, respectively. \cite{wesolowsky1973dynamic, wesolowsky1975multiperiod, galvao1992lagrangean, dias2007efficient} involve facilities' opening and closing costs in models. Moreover, \cite{ahmed2003dynamic, romauch2005dynamic, marques2013simple} present dynamic location problems under uncertain environments. \cite{nickel2019multi} is a survey of dynamic facility location problems.

In the past century, selecting an optimal retail location in a competitive environment has been much discussed. The concept of stability in the competition was first proposed in \cite{hotbllino1929stability}. \cite{reilly1931law, converse1949new} are early works that mention the retailing models about the relationship between customer consumption, the size of the facility, distance, et cetera. \cite{huff1964defining, huff1966programmed} presents further progress: Huff's formulation of the competitive function of allocating to customers; depends on the facility's attraction and distance to the customer (\cite{jiang2019solving,liang2020calibrating, jiang2021using} collect social media data to evaluate the facility's attraction to coincide with the current trend). On account of Huff's formulation, \cite{nakanishi1974parameter} adds more possible factors to extend the original formulation with the least squares approach. \cite{drezner2023competitive} is a recent paper discusses circumstances in the competitive environment to avoid cut-throat peer competition that eliminates existing facilities if a new facility is allocated. To sum up, \cite{ghosh1987location, eiselt1993competitive, ashtiani2016competitive, eiselt2019competitive} are surveys of the competitive location models.

Emerging as the times require, applications of relocatable facilities have the benefits of adapting to the times to optimize the cost and meet different needs. Many multifarious applications have been studied for relocatable facilities. \cite{cho2021spatiotemporal,bhandawat2018location,kalra2022intelligent} study location optimization problems for mobile vendors or food trucks in interurban environments to serve more customers. \cite{adler2014location} focuses on allocating police routine patrol vehicles on the city streets to improve public security. \cite{contardo2012balancing} investigates practical problems for the bike-sharing system: there is no available bike for renting and capacitated stations without an empty spot for returning at rush hour. The article provides a solution to the dynamical model to balance the inventory of unused bikes among stations with adequate inventory and short inventory. \cite{ruaboacua2020optimization} lies their research on how to temporarily locate transportable charging stations for the electric vehicle under potential constraints with optimal charging demand. \cite{glaeser2019optimal} interrelates to the distribution problem of E-commerce. Due to the high charge of door-to-door delivery, there is an innovative solution that customers can order online first, then pick up offline at a location where delivery trucks allow parking and maximize the profit of the logistics company. An analogous application is given in \cite{cao2022stall} about the unmanned market on wheels for stall economy; dis-similarly, the emphasis is on planning routes based on the reformative vehicle problem model.

\section{Modelling Process and Methods} \label{sec:model}

This section presents the development of a mixed-integer linear programming (MILP) model to address the research problem.

\subsection{Notation list}
The following are the notations used in the formulation of the model, including sets, parameters, and decision variables:

\begin{table}[H]
\centering
\caption{Notations used for mathematical modeling of DMS-p-MP}
\label{tab:groups}
\begin{adjustbox}{width=\textwidth}
\begin{tabular}{l l} 
\hline
Sets &  \\ 
\hline\hline
$T$ & Set of time periods in the planning horizon; $ t \in \{1,2,\dots,|T|\}$  \\ 
$K$ & Set of categories (groups) in the area\\
$B$ & Set of candidate locations; $ i,j \in \{1,2,\dots,|B|\}$ \\
&  where based on our problem, the candidate locations ($j$) are the same as demand nodes ($i$). \\
\hline
Parameters & \\
\hline\hline
$c_{ij}$ & Unit cost of satisfying demand of location $i$ from facility $j$ \\
$\gamma^o$ & Mobile store's opening cost in each location\\
$\gamma^c$ & Mobile store's closing cost in each location\\
$d_{it}$ & Demand of location $i$ at day $t$\\
$p$ & Available number of mobile stores in each day \\
$m_k$ & Maximum number of stores allowed in group $k$\\
$n_k$ & Minimum number of stores allowed in group $k$\\
\hline
Decision Variables & \\
\hline\hline
$x_{ijt}$ & Fraction of demand of $i$ that is supplied from $j$ at day $t$ \\ 
$y_{jt}$ & Binary variables that is 1 if a mobile store is located at $j$ at day $t$, and is 0 otherwise \\
$a_{jt}$ & Auxiliary binary variables which 1 if a store is located in $j$ \\& at day $t$ and will not be located in $j$ at day $t+1$ (,i.e., closing variable), and is 0 otherwise\\
$b_{jt}$ & Auxiliary binary variable which is 1 if a store is not located in $j$\\& at day $t$ and will be located in $j$ at day $t+1$ (,i.e., opening variable), and is 0 otherwise\\
\hline
\end{tabular}
\end{adjustbox}
\end{table}

The set of all candidate locations (buildings) in group $k$ is denoted by $F_k$. We know $B = \cup (F_k)_{k\in [K]}$, that it means each building should be at least in one group. Note that $\cap (F_k)_{k\in [K]}$ can be nonempty, meaning that one building can be a member of more than one group. Moreover, we do not need to define opening and closing costs ($a_{jt}$, $b_{jt}$) as binary variables since our minimization model and constraints force them to be 0 or 1. In order to have a feasible solution, we assume that $p$ must satisfy the following inequalities:
\begin{align*}
    \sum_{k\in K} n_k \leq p \leq \sum_{k\in K} m_k,
\end{align*}
meaning that there is a solution that satisfies the model's assumptions.

\subsection{Mathematical Formulation}
We defined all parameters and variables of our problem. The goal is to minimize the total cost of the system. We can model our problem as follows:

\begin{mini!}
{}{ \sum_{t\in T}\sum_{j\in B}\sum_{i\in B}d_{it}c_{ij}x_{ijt}+\sum_{t\in T\setminus\{|T|\}}\sum_{j\in B}(\gamma^{c}a_{jt}+\gamma^{o}b_{jt})}{\label{objective}}{\text{Optimal Cost} =}
\addConstraint{\sum_{j\in B}x_{ijt}=1} {,\label{con1}}{}{\mkern53mu \forall i\in B, \, t\in T}
\addConstraint{\sum_{j\in B}y_{jt} =p} {,\label{con2}}{}{\mkern53mu \forall t\in T}
\addConstraint{x_{ijt} - y_{jt}\leq 0} {,\label{con3}}{}{\mkern53mu \forall i,j\in B, \, t\in T}
\addConstraint{y_{jt}-y_{j(t+1)} -a_{jt}\leq 0} {,\label{con4}}{}{\mkern53mu \forall j\in B,\, t\in T\setminus\{|T|\}}
\addConstraint{y_{j(t+1)}-y_{jt}-b_{jt}\leq 0} {,\label{con5}}{}{\mkern53mu \forall j\in B,\, t\in T\setminus\{|T|\}}
\addConstraint{\sum_{j \in F_{k}}y_{jt}\leq m_k} {,\label{con6}}{}{\mkern53mu \forall k\in K,\, t\in T}
\addConstraint{\sum_{j \in F_{k}}y_{jt}\geq n_k} {,\label{con7}}{}{\mkern53mu \forall k\in K,\, t\in T}
\addConstraint{y_{jt}\in\{0,\,1\}} {,\label{con8}}{}{\mkern53mu \forall j \in B,\, t\in T}
\addConstraint{a_{jt},\, b_{jt}\geq 0} {,\label{con9}}{}{\mkern53mu \forall j\in B, \, t\in T}
\addConstraint{x_{ijt}\geq 0} {,\label{con10}}{}{\mkern53mu \forall i,j\in B, \, t\in T.}
\end{mini!}
where constraint (\ref{con1}) ensures the demand of all locations over the time horizon must be satisfied, constraint (\ref{con2}) ensures the number of mobile stores allocated in each day should be $p$, constraint (\ref{con3}) shows the relation between continuous and binary variables, i.e., if $x_{ijt}>0$, then $y_{jt}$ should be 1. Constraints (\ref{con4}) and (\ref{con5}) ensure that we consider the closing and opening costs, respectively. Constraints (\ref{con6}) and (\ref{con7}) ensure that we are following the policies regarding the limits on the number of stores in each group. Finally, constraint (\ref{con8}) is the binary constraint, and (\ref{con9}), and (\ref{con10}) are the non-negativity constraints. We assume that at the beginning of day 1, $d_{it}$ is stochastic for all $i\in B$ and all $t\in T$. 

\subsection{Solution Approach}
The problem at hand is a mixed-integer mathematical programming model that involves uncertainty. We will outline a state-of-the-art methodology to tackle these characteristics in the following manner:

\subsubsection{Robust Optimization}
In this section, an overview of the robust optimization approach proposed by \cite{bertsimas2004price} is presented. To do so, the following linear programming model is considered:

\begin{flalign} \nonumber 
Min_{} \quad &\sum_{j}^{} c_{j}x_{j}\\ \nonumber 
s.t. \quad &\sum_{j}^{} \Tilde{a}_{ij}x_{j} \leq b_{i}; \qquad\forall i\\
&x_{j}\geq 0; \qquad\forall j
\end{flalign}

where the technological coefficients $\Tilde{a}_{ij}$ are assumed to be uncertain. In other words, each coefficient $\Tilde{a}_{ij}$ is regarded as an independent, symmetric, and bounded parameter, which can take values in $[a_{ij} - \hat{a}_{ij},a_{ij} + \hat{a}_{ij}]$, i.e. $\Tilde{a}_{ij} \in [a_{ij} - \hat{a}_{ij},a_{ij} + \hat{a}_{ij}]$. In this definition, $a_{ij}$ and $\hat{a}_{ij}$ denote the nominal value and the maximum deviation from the nominal value, respectively. Associated with each row \textit{i} in problem (1) is $J_{i}$, which is defined as the set of all coefficients in row \textit{i} that are subject to uncertainty. Furthermore, a scaled deviation $\eta_{ij} \in [-1,1]$ is defined for each uncertain coefficient $\Tilde{a}_{ij}$ as $\eta_{ij} = \dfrac{\Tilde{a}_{ij} - {a}_{ij}}{\hat{a}_{ij}}$ that represents the scaled perturbation of $\Tilde{a}_{ij}$ from its nominal value $a_{ij}$.

Bertsimas and Sim (2004) also introduced a parameter $\Gamma_{i} \in [0,|J_{i}|]$ as the budget of uncertainty for each constraint \textit{i}, where $|J_{i}|$ denotes the number of elements of set $J_{i}$ \cite{bertsimas2004price}. In fact, $\Gamma_{i}$ is the maximum number of parameters that can really deviate from their nominal values for each constraint \textit{i}. The parameter $\Gamma_{i}$ that bounds the total scaled deviation of uncertain parameters as $\sum_{j \in J_{i}} |\eta_{ij}| \leq \Gamma_{i}$ adjusts the robustness of the proposed method against the level of solution conservatism. In particular, $\Gamma_{i}=0$ represents the nominal or deterministic formulation, whereas $\Gamma_{i}=|\eta_{ij}|$ relates to the worst-case formulation in which all uncertain parameters are fixed at their worst-case values from the uncertainty set. However, the decision maker can make a trade-off between the protection level of constraint \textit{i} and the degree of conservatism of the solution if $\Gamma_{i} \in (0,|J_{i}|)$. Therefore, the budget of uncertainty $\Gamma_{i}$ that is an input to the robust optimization model can specify how risk-averse the decision-maker is. 

Bertsimas and Sim (2004) proposed a nonlinear programming model as follows, which is equivalent to the uncertain model (1):
\begin{flalign} \nonumber 
Min_{} \quad &\sum_{j}^{} c_{j}x_{j}\\ \nonumber 
s.t. \quad &\sum_{j}^{} a_{ij}x_{j} + \underset{\Omega}{max}\{\sum_{j \in S_{i}}^{} \hat{a}_{ij}x_{j} + (\Gamma_{i} - \left \lfloor \Gamma_{i} \right \rfloor)\hat{a}_{it_{i}}x_{j}\} \leq b_{i}; \qquad\forall i\\
&x_{j}\geq 0; \qquad\forall j 
\end{flalign}
where $\Omega=\{S_{i}\cup\{t_{i}\}|S_{i}\subseteq J_{i}, S_{i}=\left \lfloor \Gamma_{i} \right \rfloor, t_{i} \in J_{i} \setminus S_{i} \}$ is defined as the uncertainty set. For a given optimal solution $x^{*}$ of the problem (2), Bertsimas and Sim (2004) demonstrated that the protection function for constraint \textit{i} against uncertainty, which is $\beta_{i}(x^{*},\Gamma_{i})=\underset{\Omega}{max}\{\sum_{j \in S_{i}}^{} \hat{a}_{ij}x_{j} + (\Gamma_{i} - \left \lfloor \Gamma_{i} \right \rfloor)\hat{a}_{it_{i}}x_{j}\}$ can be formulated as the following linear programming problem:
\begin{flalign} \nonumber 
\beta_{i}(x^{*},\Gamma_{i}) = Max_{} \quad &\sum_{j \in J_{i}} \hat{a}_{ij}|x_{j}^{*}|\eta_{ij}\\ \nonumber 
s.t. \quad &\sum_{j \in J_{i}}^{} \eta_{ij} \leq \Gamma_{i}; \qquad\forall i\\
&0 \leq \eta_{ij}\leq 1; \qquad\forall i,j 
\end{flalign}
According to the theory of strong duality, since problem (3) is always feasible and bounded for all $\Gamma_{i} \in [0,|J_{i}|]$, its dual problem is feasible and bounded as well. Therefore, replacing the dual problem of the problem (3) into (2), Bertsimas and Sim (2004) derived the robust formulation of the uncertain linear programming problem (1) as follows: 
\begin{flalign} \nonumber 
Min_{} \quad &\sum_{j}^{} c_{j}x_{j}\\ \nonumber 
s.t. \quad &\sum_{j}^{} a_{ij}x_{j} + \beta_{i}\Gamma_{i} + \sum_{j \in J_{i}}\mu_{ij} \leq b_{i}; \qquad\forall i\\ \nonumber
& \beta_{i} + \mu_{ij}\geq \hat{a}_{ij}x_{j}; \qquad\forall i, j \in J_{i}\\ \nonumber 
& \mu_{ij} \geq 0; \qquad\forall i, j \in J_{i}\\ \nonumber
&\beta_{i} \geq 0; \qquad\forall i\\
&x_{j}\geq 0; \qquad\forall j 
\end{flalign}
where $\beta_{i}$ and $\mu_{ij}$ are dual variables associated with the first and second constraints in programming problem (3), respectively.

if the number of uncertain coefficients in constraint \textit{i} that perturb from their respective nominal values is less than or equal to $\Gamma_{i}$, then the optimal solution from the robust problem (4) will remain always feasible. However, if more than $\Gamma_{i}$ coefficients deviate from their nominal values, then the probability of violating constraint \textit{i} for an optimal solution $x_{j}^{*}$ is calculated as follows:

\begin{flalign} 
Pr(\sum_{j} \Tilde{a}_{ij}x_{j}^{*} < b_{i}) \leq 1-\varphi(\dfrac{\Gamma_{i} - 1}{\sqrt{|J_{i}|}}) 
\end{flalign}
where $\varphi(.)$ is the cumulative distribution function of a standard normal random variable. 

The robust optimization equivalent of the proposed problem (\ref{objective}) is:
\begin{subequations}
\begin{align}
    \min_{x,y,a,b} \max_{d\in \mathcal{U}}\ & {d}{c}^1{x} +{c}^2(a+b)\\
    \text{s.t. } \ &(1b)-(1k)
\end{align}
\end{subequations}
or equivalently:
\begin{subequations}
\begin{align}
    \min_{x,y,a,b} & \theta +{c}^2(a+b)\\
    \text{s.t. } \ &(1b)-(1k)\\
    & \theta \geq \max_{d\in \mathcal{U}} {d}{c}^1{x}
\end{align}
\end{subequations}

where $\mathcal{U}$ is the uncertainty set for demand, and $a$, $b$ are decision vectors for $a_{jt}$, $b_{jt}$ respectively. Depending on the problem setting, one can use the box, budget, or conic uncertainty set. Among these, the budget uncertainty set is widely used, both due to its intuitive interpretation and tractability. The budget uncertainty set for the demand is:

$$\mathcal{U}= \{ {d} \in \mathbb{R}^{|\mathcal{D}|\times |\mathcal{D}|}~ | ~ {d} = \overline{{d}} + {\gamma}{\xi}, {\xi}||_{\infty} \leq 1, \sum_{i,j\in \mathcal{D}} \xi_{it}\leq \Gamma \} $$ 
 where $\Gamma$ is the uncertainty budget. Using this definition we have:
\begin{subequations}
\begin{align}
    \min_{x,y,a,b} & \theta +{c}^2(a+b)\\
    \text{s.t. } \ &(1b)-(1k)\\
    &\max_{{\xi}_\infty \leq 1, {\xi}_1 \leq \Gamma} {\xi}{c}^1{x}\leq  \theta - {\bar{d}}{c}^1{x}
\end{align}
\end{subequations}
Let's focus on the left-hand-side of equation (3c):
\begin{subequations}
    \begin{align}
        &\max_{{\xi}_\infty \leq 1, {\xi}||_1 \leq \Gamma}  \gamma \xi c^1 x \to \max_{{\xi}} \gamma \xi c^1 x, \quad  s.t.: \ \sum_{i,t} |\xi_{it}|\leq \Gamma, -1 \leq \xi_{it}\leq 1 
    \end{align}
\end{subequations}
Let $a_{it} = |\xi_{it}|$, then we have:

\begin{subequations}
    \begin{align}
       \max \ &  \gamma \xi c^1 x\\
       \text{s.t.} \ & \sum_{i,t} a_{it}\leq \Gamma & [\pi^1]\\
       & \xi_{it}\leq 1 & [\pi^2_{it}]\\
       & -\xi_{it}\leq 1 &[\pi^3_{it}]\\
       & \xi_{it} - a_{it} \leq 0&[\pi^4_{it}]\\
       & -\xi_{it} -  a_{it} \leq 0& [\pi^5_{it}]\\
       & \xi_{it}\in \mathbb{R}, a_{it}\in \mathbb{R}^+
    \end{align}
\end{subequations}

Since the model is a minimization problem, we need to obtain the dual problem whose objective is also minimization:
The dual problem is:
\begin{subequations}
    \begin{align}
        \min_{\boldsymbol{\pi}} \ & \Gamma \pi^1 +\sum_{i,t\in \mathcal{D}} (\pi^2_{it}+\pi^3_{it})\\
        \text{s.t.} \quad \ & \pi^1 \geq \pi^4_{it} + \pi^5_{it} \quad \forall i,t\in \mathcal{D}&\\
        & \pi \in \mathbb{R}^+
    \end{align}
\end{subequations}

Finally the DMS-p-MP reformulation of the problem becomes as:
\begin{subequations}
\begin{align}
    \min_{x,y,a,b} & \theta +{c}^2(a+b)\\
    \text{s.t. } \ &(1b)-(1k)\\
    & \Gamma \pi^1 +\pi^2+\pi^3 + \Bar{d}c^1 x \leq \theta\\
    & \pi^1 \geq \pi^4_{it} + \pi^5_{it} \quad \forall i,t\in \mathcal{D} &\\
        & \pi^2_{it}-\pi^3_{it}+\pi^4_{it} - \pi^5_{it} = \gamma_{it} c^1_{it} x_{ijt} \quad \forall i,t\in \mathcal{D}, \forall j \in B &\\
        & \pi^\kappa_{it} \geq 0 \quad \forall i \in B, t \in T, \kappa \in \{1,2,3,4,5\}
\end{align}
\end{subequations}

where $\pi^\kappa$, $x$, and $d$  is the vector of $\pi^\kappa_{it}$'s, $d_{it}$, and $x_{ijt}$'s respectively.

\subsubsection{Lagrangian relaxation}
In discrete location theory, one of the basic models is the p-median problem and it is an NP-hard problem, as with most location problems \cite{mladenovic2007p}. The MILP model presented in section \ref{sec:model} is typically solved in practice with numerous residence areas (demand points), potential facility locations, and time horizons. As a result of the large-scale property of the developed model, it presents substantial computational difficulty when applied to real problems, which cannot be addressed by commercial optimization software, such as Cplex, Xpress or Gurobi. For supply chain optimization problems involving such computational complexity, Lagrangian relaxation (LR) is widely used \citep{diabat2013lagrangian, duong2018mixed, rafie2018modelling, kheirabadi2019mixed, hamdan2020robust} in the literature.

Therefore, this section develops an LR approach for solving the presented MILP problem. The LR method is an iterative algorithm that provides the upper and lower bounds of the optimal objective value as well as the estimation of the optimality gap of the feasible established solution in each iteration \citep{daskin1997network}.

The LR method used in this paper includes the general steps as follows: 1) Relax one of the constraints by multiplying it by a Lagrange multiplier and bringing the constraint into the objective function, 2) Solve the model to find the optimal values of the relaxed problem,
3) Find the feasible solution to the original problem by using the resulting decision variables found in step 2, 4) Compute the lower bound using the solution obtained from the relaxed problem in step 2, 5) Use the subgradient optimization method to modify the Lagrange multiplier assigned to the violated constraint, return to step 2 after finding the new multiplier(s) for the Lagrange variable. The algorithm terminates whenever the lower bound is close enough to the upper bound.

\paragraph{Step 1. Solving the Relaxed Problem}
In order to make the problem easier to solve, the constraints in this study are relaxed, even if this relaxation may lead to infeasibility \citep{daskin1997network}. Specifically, constraints (1b), (1c) are relaxed, resulting in the following Lagrangian dual problem:

\begin{subequations}
\begin{align}
\nonumber    \min_{x,y,a,b} & \theta +{c}^2(a+b)+\\& \nonumber
 \sum_{j \in B} \sum_{i \in B} \sum_{t \in T} \lambda_{bt}^{1}(x_{ijt}-1)+\\& 
 \sum_{j \in B} \sum_{t \in T} \lambda_{t}^{2}(y_{jt}-p)\\ 
    \text{s.t. } \ &\text{(1d) - (1k), (13c) - (13f)}
\end{align}
\end{subequations}

The optimal value of the objective function in the Lagrangian dual problem above, which is defined using non-negative Lagrange multipliers ($\lambda_{bt}^{1}$, $\lambda_{t}^{2}$), serves as a lower bound for the mixed-integer linear programming problem with constraints numbered (14a) to (14g).

\paragraph{Step 2. Finding a Feasible Solution and an Upper Bound}
In most situations, the solution to the Lagrangian dual problem is not feasible due to the relaxation of constraints (1b) and (1c). However, it is possible to obtain a feasible solution that provides an upper bound on the single objective MILP model by solving this model and setting the decision variables $x_{ijt}$, $y_{jt}$, $a_{jt}$ and $b_{jt}$ to the optimal values obtained from solving the Lagrangian dual problem.

\paragraph{Step 3. Finding a Lower Bound, and Updating the Lagrange Multipliers}
During each iteration of the Lagrangian procedure, the Lagrangian multipliers $\lambda_{bt}^{1}$, $\lambda_{t}^{2}$ are updated and new lower and upper bounds are subsequently derived. There are various methods in the literature, such as cutting planes \citep{kelley1960cutting}, sub-gradients \citep{daskin1997network}, and bundling \citep{borghetti2003lagrangian}, that can be used for this purpose. In this paper, the sub-gradient approach is employed to update the Lagrangian multipliers because it is a widely recognized and commonly used method. According to the sub-gradient procedure \citep{daskin1997network}, the Lagrange multipliers at the $(n + 1)$ iteration are calculated as follows:

 \begin{flalign}
        \label{eq:LRupdate1}  & (\lambda_{bt}^{1})^{n+1} = \max \left \{0,(\lambda_{bt}^{1})^{n} - \tau_{1}^{n}(x_{ijt}-1) \right\} \\
        \label{eq:LRupdate2}  & (\lambda_{t}^{2})^{n+1} = \max \left \{0,(\lambda_{t}^{2})^{n} - \tau_{2}^{n}(y_{jt}-p) \right\}  && 
\end{flalign}

The step sizes $\tau_{1}^{n}$, and $\tau_{2}^{n}$ in the algorithm are defined as follows:
 \begin{flalign}
        \label{eq:LRstepSize1}  & \tau_{1}^{n} = \frac{\alpha^{n}(UB-LB^{n})}{\sum_{j \in B} \sum_{i \in B} \sum_{t \in T} \lambda_{bt}^{1}(x_{ijt}-1)^{2}} \\
        \label{eq:LRstepSize2}  & \tau_{2}^{n} = \frac{\alpha^{n}(UB-LB^{n})}{\sum_{j \in B} \sum_{t \in T} \lambda_{t}^{2}(y_{jt}-p)^{2}}  && 
\end{flalign}

The term $\alpha$ is simply a constant that will be changed during each iteration of the algorithm as described above. $\alpha$ is initialized to 2. If the lower bound, $LB$, has not increased in an iteration, then its value will be halved. Additionally, let $UB$ be the best upper bound (the one with the smallest value) that has been discovered thus far, and $LB^n$ is the lower bound obtained at iteration $n$.

\paragraph{Step 4. Termination Criteria}
The algorithm will stop when one of the following conditions is met \citep{daskin1997network}:

\begin{itemize}
    \itemsep=0pt
      \item A predetermined number of iterations have been completed.
      \item The lower bound is equal to the upper bound ($UB = LB^n$) or is close enough to the upper bound ($UB-LB^n < 0.1$).
      \item The value of $\alpha$ becomes small
\end{itemize}

\section{Numerical Experiment}\label{sec:example}
\subsection{Case description}

\begin{figure}[h!]
 \centering 
 \includegraphics[width=12cm]{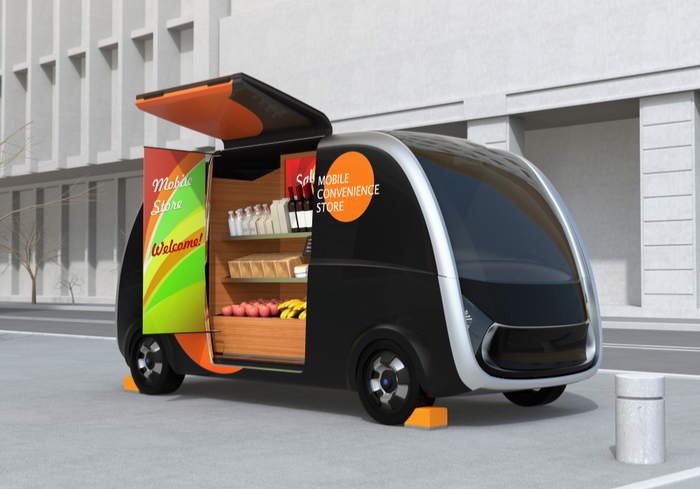}
 \caption{Sample Picture of Mobile Grocery Store from \cite{pymnts.com_2018}}
\label{robo}
\end{figure}

In this section, we discuss a case in which a company wants to place some mobile grocery stores on the campus of the University of Waterloo. These self-service mobile stores can serve all students and staff. Figure \ref{robo} shows a sample picture of these stores. Specifically, the company wants to propose a plan to specify the optimal locations of stores needed on campus dynamically (every day) based on the demand variation on campus. They update their plan every month, i.e., they decide on all days at the beginning of each month. In this research, we aim to find the optimal locations of stores over the time horizon (one month). Based on Section~\ref{sec:model}, $|T|=28$, $|K|=6$, and $|B|=91$. 

\begin{figure}[H]
 \centering 
 \frame{\includegraphics[width=12cm]{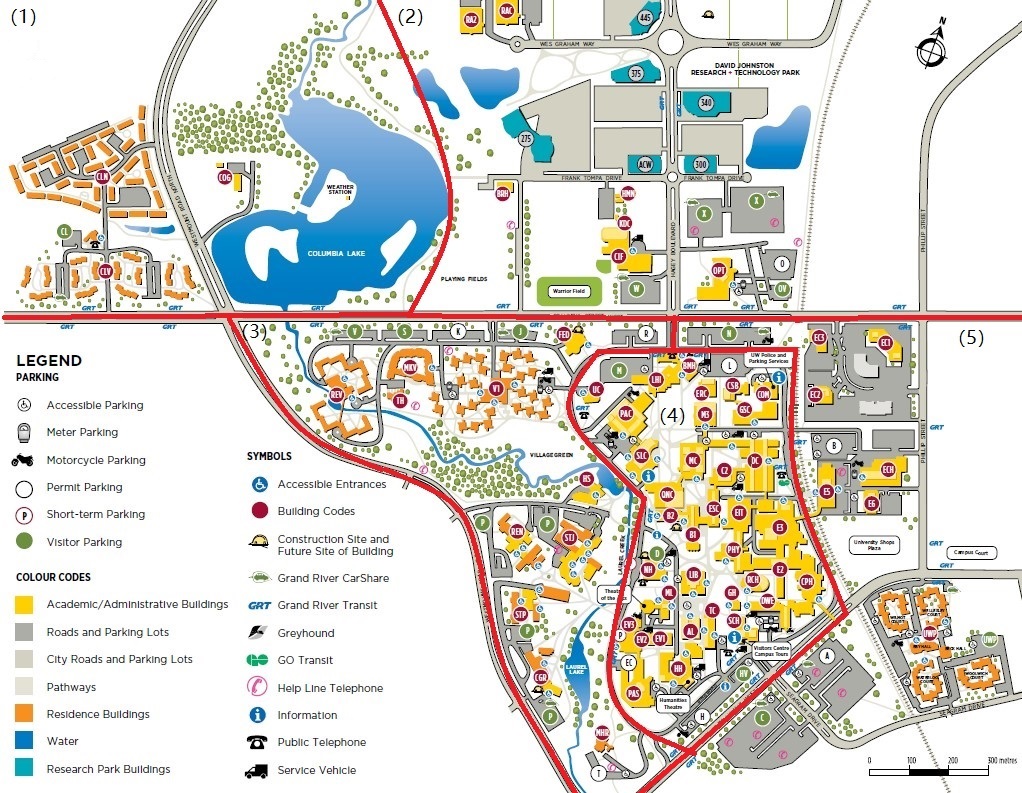}}

 \caption{University of Waterloo - Campus Map}
\label{uwmap}
\end{figure}

Figure \ref{uwmap} shows the University of Waterloo campus map. All buildings on campus are categorized into five groups: Service and Administrative Buildings, Academic Buildings, Residence Buildings, Research Park Buildings, and others. The University of Waterloo asked the company to follow some rules regarding each group. Specifically, there are some limitations on the number of stores needed for each group for each day. More details and the list of buildings included in each group will be provided in the following sections.

Each building has a specific demand on each day. Obviously, it is not (financially) feasible to place one store in each building. The company has only $p$ stores that in each period (day) they want to place all of them on campus. We consider a \emph{demand cost}, for students and staff in a building that does not have a store and they need to go to other buildings. The goal is to place the stores to satisfy the demand of all buildings while minimizing the costs. There are two more costs in our problem: \emph{opening cost} and \emph{closing cost}. Although the stores are mobile, we need to consider an opening (closing) cost if we want to open (close) a store at the beginning of each day. In fact, we are capturing all transportation costs needed to change a store's location. 

The demand of each building may change day-to-day because of different reasons such as weekends, events, etc. Then, it is an excellent motivation for day-to-day decision-making following demand variation. In this research, we assume that our demand prediction is accurate. Then, we can finalize all decisions for each day of the month, once at the beginning of the month, i.e., our demand is deterministic.

We can model this problem as a \emph{modified} multi-period $p$-median problem. Generally, our problem has two modifications compared to the simple multi-period $p$-median problem. First, we are considering opening and closing costs over the time horizon. Second, there are some extra constraints raised from University of Waterloo policies.

\subsection{Data Collection}
We divide the University of Waterloo's (UW's) all facilities according to their functionality into six different segments: Academic Buildings, Parking Spots, Student Residence Buildings, Research Park Buildings, Athletic Buildings, and University Plaza respectively. The workspace is based on the UW's official map includes 91 facilities, and we also make an assumption that the population in each facility within a segment is all equivalent. Segments indicate in the Table \ref{segments} are as follows,

\begin{table}[H]
\centering
\caption{The Functional Classification for All Facilities in the Workspace}
\label{segments}
\begin{tabular}{l p{8cm}} 
 \hline
 Segment  & Buildings Included in the Segment (Building Id)\\ 
 \hline\hline
Academic Buildings & COG (1), COM (2), CPH (3), RA2 (4), M3 (5), ML (6), RAC (7), GSC (8), GH (9), BRH (10), SLC (11), OWE (12), HMN (13), MC (14), TC (15), KOC (16), C2 (17), EV3 (18), OPT (19), DC (20), SCH (21), FED (22), QNC (23), AL (24), HS (25), B2 (26), EV2 (27), REN (28), ESC (29), EV1 (30), STJ (31), EIT (32), HH (33), STP (34), Bl (35), PAS (36), CGR (37), E3 (38), EC3 (39), BMH (40), PHY (41), EC1 (42), LHI (43), NHI (44), EC2 (45), UC (46), LIB (47), ECH (48), ERC (49), E2 (50), ES (51), CSB (52), RCH (53), E6 (54).\\
\hline
Parking	& Parking CL (55), Parking A (56), Parking X (57), Parking C (58), Parking W (59), Parking OV (60), Parking V (61), Parking S (62), Parking K (63), Parking J (64), Parking R (65), Parking P (66), Parking T (67), Parking M (68), Parking L (69), Parking D (70), Parking EC (71), Parking HV (72), Parking N (73), Parking UWP (74).\\
\hline
Residence Buildings & CLN (75), CLV (76), MKV (77), V1 (78), REV (79), TH (80), MHR (81), UWP (82).\\
\hline
Research Park Buildings	& 445 (83), 375 (84), 340 (85), 275 (86), ACW (87), 300 (88).	\\
\hline
Athletic Buildings & CLF (89), PAC (90).\\
\hline
University Plaza & Plaza (91).\\
\hline
\end{tabular}

\end{table}

\begin{figure}[H]
  \label{tbl:excel-table}
  \begin{center}
      \frame{\includegraphics[width=14cm]{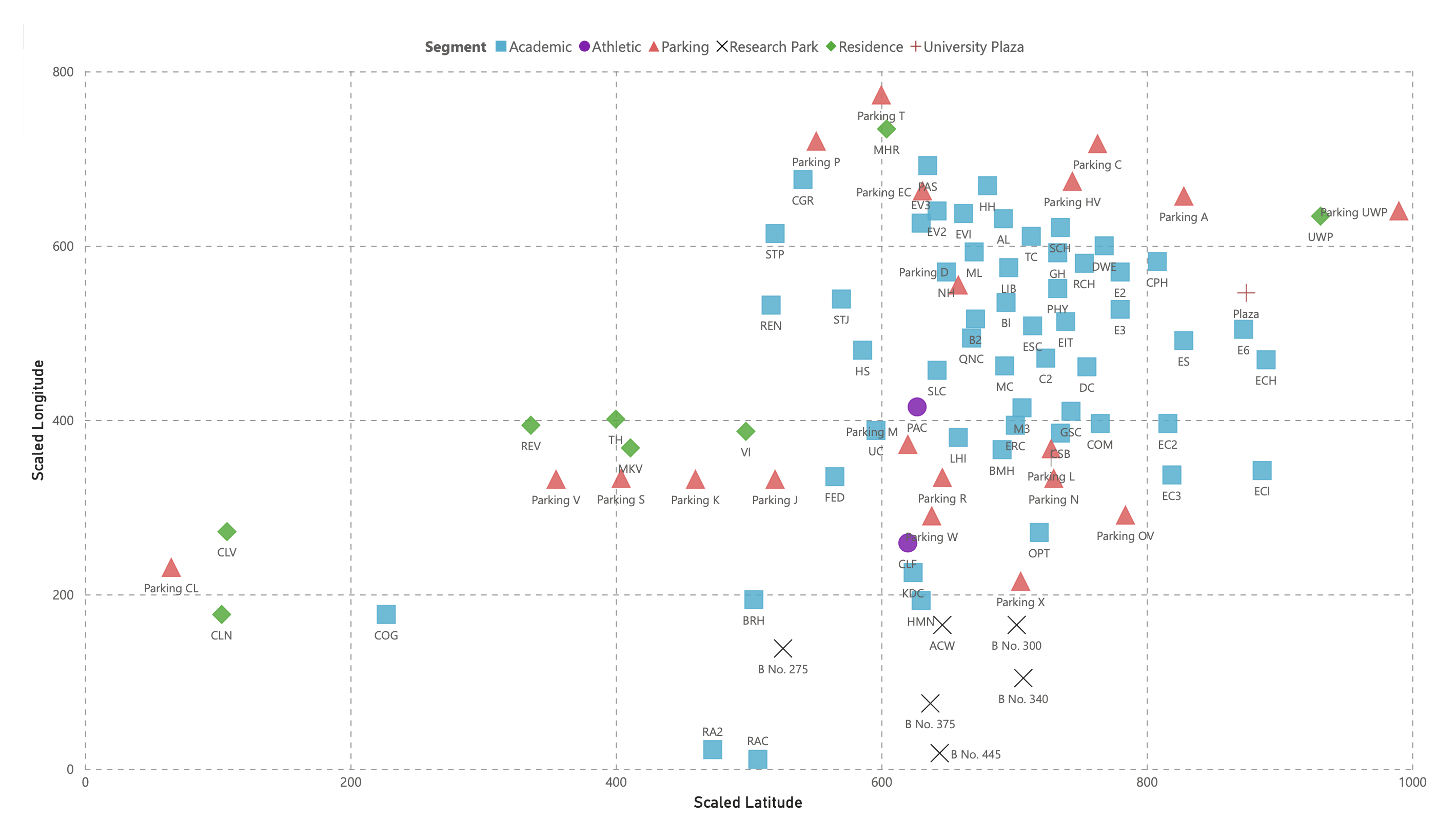}}
  \end{center}
   \caption{Coordinates for All Facilities on the Workspace}
   \label{Coordinates}
\end{figure}

We use Python's OpenCV toolbox by \cite{bradski2000opencv} to determine facilities' specific locations by fixing pixel's coordinates in the workspace. Students and staff allows going through buildings, it is likely for us to use Euclidean distance measurement
$$l_{2}=\sqrt{(x-a)^{2}+(y-b)^{2}}$$ to find the distance between the two pairs in the data set is as follows on the Figure~\ref{Coordinates}.

From the UW's website, we obtained the rough number of students and staff for each faculty, the rough number of students and staff for each university college, and the capacities of each residence. A general assumption is that 10\% of university students will use one of two gyms (CIF and PAC). Tables \ref{Segment1}, \ref{Segment2}, \ref{Segment3}, \ref{Segment4}, \ref{Segment5}, \ref{Segment6} below are the population of all classifications for each segment.

\begin{table}[H]
\caption{Statistical Results for Academic Buildings}
\label{Segment1}
\label{tab:groups}
    \centering
\begin{tabular}{c c c} 
 \hline
 Segment (1) & Population &  \\ 
 \hline\hline
 Engineering Faculty &11000 &  \\ 
 \hline
 Mathematics Faculty&9260 & \\
 \hline
 Science Faculty & 6000&\\
 \hline
 Health Faculty & 3643&\\
 \hline
 Arts Faculty &3000&\\
 \hline
 Environment Faculty &3000&\\
 \hline
 Others & 1200 &\\
 \hline
 \hline
 Total population&37103&\\
 \hline
 Total facilities &54&\\
 \hline
 Average per facility &687&\\
 \hline
\end{tabular}
\end{table}

\begin{table}[H]
\caption{Statistical Results for Parking Spots}
\label{Segment2}
\label{tab:groups}
    \centering
\begin{tabular}{c c c} 
 \hline
 Segment (2) & Population &  \\ 
 \hline\hline
 Campus parking &4000 &  \\ 
 \hline
 \hline
 Total population&4000&\\
 \hline
 Total facilities &20&\\
 \hline
 Average per facility &200&\\
 \hline
\end{tabular}
\end{table}

\begin{table}[H]
\caption{Statistical Results for Student Residence Buildings}
\label{Segment3}
\label{tab:groups}
    \centering
\begin{tabular}{c c c} 
 \hline
 Segment (3) & Population &  \\ 
 \hline\hline
 Columbia Lake Village - North &404 &  \\ 
 \hline
 Columbia Lake Village - South &400 & \\
 \hline
  William Lyon Mackenzie King Village& 320 & \\
 \hline
 Student Village 1&1381&\\
 \hline
  Ron Eydt Village &960 &  \\ 
 \hline
 Tutors’ Houses &100&\\
 \hline
 Minota Hagey Residence &70&\\
 \hline
 University of Waterloo Place &1650&\\
 \hline
 Others &200&\\
 \hline
 \hline
 Total population&5285&\\
 \hline
 Total facilities  &8&\\
 \hline
 Average per facility &660&\\
 \hline
\end{tabular}
\end{table}

\begin{table}[H]
\caption{Statistical Results for Research Parking Buildings}
\label{Segment4}
\label{tab:groups}
    \centering
\begin{tabular}{c c c} 
 \hline
 Segment (4) & Population &  \\ 
 \hline\hline
 David Johnson Research Park &4000 &  \\
 \hline
 Others &100&\\
 \hline
 \hline
 Total population&4100&\\
 \hline
 Total facilities &6&\\
 \hline
 Average per facility &683&\\
 \hline
\end{tabular}
\end{table}

\begin{table}[H]
\caption{Statistical Results for Athletic Buildings}
\label{Segment5}
\label{tab:groups}
    \centering
\begin{tabular}{c c c} 
 \hline
 Segment (5) & Population &  \\ 
 \hline\hline
 Columbia Icefield &2100 &  \\ 
 \hline
 Physical Activities Complex &2100& \\
 \hline
 Others & 100 &\\
 \hline
 \hline
 Total population&4300&\\
 \hline
 Total facilities &2&\\
 \hline
 Average per facility &2150&\\
 \hline
\end{tabular}

\end{table}

\begin{table}[H]
\caption{Statistical Results for University Plaza}
\label{Segment6}
\label{tab:groups}
    \centering
\begin{tabular}{c c c} 
 \hline
 Segment (6) & Population &  \\ 
 \hline\hline
 University Shops Plaza &3200 &  \\ 
 \hline
 \hline
 Total population&3200&\\
 \hline
 Total facilities &1&\\
 \hline
 Average per facility &3200&\\
 \hline
\end{tabular}
\end{table}
Given in the previous section that the time horizon $T=28$ (approximately equivalent to a month), it can also separate into four periods (weeks). Based on the functionality of each segment, with time-varying from Monday to Sunday, we may empirically define the facility in the different segments and the different days in utilization rate $U_{kt}$. For instance, the utilization rates for the academic buildings are 100 out of 100, but there are 30 out of 100 during the weekend. The complete estimated utilization rates for different functional buildings for each day over a week are as the following Table \ref{utilization}:

\begin{table}[H]

\caption{The Estimated Utilization Rate for Different Functional Buildings for Each Day over a Week}
\label{utilization}
\begin{adjustbox}{width=\textwidth}
\begin{tabular}{llllllll}
\hline
Functionality & Monday & Tuesday & Wednesday & Thursday & Friday & Saturday & Sunday \\
\hline \hline
Academic Buildings & 100 & 90 & 90 & 80 & 90 & 30 & 30 \\
Parking Spots & 100 & 100 & 100 & 100 & 100 & 20 & 20 \\
Residence Buildings & 50 & 50 & 50 & 50 & 60 & 100 & 100 \\
Research Park Buildings & 100 & 90 & 90 & 80 & 90 & 10 & 10 \\
Athletic Buildings & 50 & 50 & 50 & 60 & 50 & 100 & 100 \\
University Plaza & 100 & 100 & 70 & 80 & 90 & 50 & 50\\
\hline
\end{tabular}
\end{adjustbox}
\end{table}

After that, the demand $d_{kt}$ of each building in segment $k$ on day $t$, where $t\in \{1,\dots,7\}$, in a week is able to define as, 
$$d_{kt}=U_{kt}\times \frac{\text{Total Population}(k)}{\text{Total Facilities}(k)}.$$

As for having a unified standard, we calculate the lower bound of opening facilities in constraint (\ref{con7}) by taking a floor of 70\% of the number of facilities in segment $k$ over the total facilities on the workspace, 
 $$n_{k}=\left\lfloor p \times \frac{\text{the Number of Facilities in Segment (k)}}{91}\times 70\%\right\rfloor.$$
 For example, $n_{1}=\left\lfloor 18 \times \frac{54}{91}\times 70\%\right\rfloor=\left\lfloor{7.47}\right\rfloor=7 .$ Similarly, the upper bound of opening facilities in constraint (\ref{con6}) by taking a ceiling of 130\% of the number of facilities in segment $k$ over the total facilities on the workspace,
 $$m_{k}=\left\lceil p\times\frac{\text{the Number of Facilities in Segment (k)}}{91}\times 130\%\right\rceil.$$
 
 The following Table \ref{bounds} is the maximum (minimum) number of Robomarts are allowed in segment k,

\begin{table}[H]
\caption{Bounds on the number of Robomarts}
\label{bounds}
\label{tab:groups}
    \centering
\begin{tabular}{c c c} 
 \hline
 Functionality & Min ($n_{k}$) & Max ($m_{k}$) \\ 
 \hline\hline
Academic Buildings & 7 & 14 \\ 
 \hline
 Parking Lots & 2 & 6 \\
 \hline
 Residence Buildings & 1 & 3 \\
 \hline
 Research Park Buildings & 0 & 2 \\
 \hline
 Athletic Buildings & 0 & 1 \\
 \hline
 University Plaza & 0 & 1\\
 \hline
\end{tabular}
\end{table}

\subsection{Computational Experiments}
In this section, we use all parameters discussed in the Data Collection section, to solve our problem. All numerical experiments have been run on an Apple M1 processor, limited to 16 GB of RAM. Gurobi has access to 8 physical cores, 8 logical processors, using up to 8 threads. This model is MILP and has $M (= 239512)$ variables and $N (= 239694)$ constraints (excluding sign constraints). We use Gurobi version 9.5.1 by \cite{bixby2007gurobi} to solve this optimization model. Since default settings in Gurobi generally work well, we are keeping all settings as default. Specifically, we use Gurobi's API embedded in Python. 

MILP models are generally solved using a linear-programming based branch-and-bound algorithm. The Gurobi provides advanced implementations of
the latest MILP algorithms including deterministic parallel, nontraditional search, heuristics, solution improvement, cutting planes, and symmetry breaking.

Based on section 5, the demand is weekly periodic. Then, we expect the model to make the same decision over the weeks. Table \ref{tab:results} shows what buildings are open during the planning horizon. We just show the days that we change our decisions. For instance, the buildings that are open between day 1 and day 6 are the same.

Table \ref{tab:results} shows that we only open new buildings and close the current buildings over the weekend. We again change our decision on weekdays. It is expected because the utilization rates of some of our segments are significantly different during the weekend.

Specifically, we close buildings 12, 55, and 83 (that base on Table \ref{Coordinates} are Conrad Grebel university college, Parking lots, Research building 375, respectively) and open buildings 63, 77, and 81 (that are Parking P, Student Village, University of Waterloo Place respectively). The interesting point is that all buildings 63, 77, and 81 are  around students' residences. Students are mostly in the residence area instead of the academic campus during the weekends. Then, it is worth closing some stores and opening new ones in students' residences.

\begin{table}[H]
\caption{Results - Optimal Objective Value = 43175.7, CPU Time = 12.99 Seconds}
\label{tab:results}
    \centering
\begin{tabular}{c p{5.5cm}} 
 \hline
 t & Open Buildings \\ 
 \hline\hline
1 (Monday) & 0, 2, 3, 5, 10, 12, 17, 27, 34, 40, 44, 48, 54, 55, 76, 83, 89, 90  \\ 
 \hline
 6 (Saturday) & 0, 2, 3, 5, 10, 17, 27, 34, 40, 44, 48, 54, 63, 76, 77, 81, 89, 90 \\ 
 \hline
 8 (Monday) & 0, 2, 3, 5, 10, 12, 17, 27, 34, 40, 44, 48, 54, 55, 76, 83, 89, 90 \\ 
 \hline
 13 (Saturday) & 0, 2, 3, 5, 10, 17, 27, 34, 40, 44, 48, 54, 63, 76, 77, 81, 89, 90 \\ 
 \hline
 15 (Monday) & 0, 2, 3, 5, 10, 12, 17, 27, 34, 40, 44, 48, 54, 55, 76, 83, 89, 90 \\ 
 \hline
 20 (Saturday) & 0, 2, 3, 5, 10, 17, 27, 34, 40, 44, 48, 54, 63, 76, 77, 81, 89, 90 \\ 
  \hline
   22 (Monday) & 0, 2, 3, 5, 10, 12, 17, 27, 34, 40, 44, 48, 54, 55, 76, 83, 89, 90 \\ 
  \hline
   27 (Saturday) & 0, 2, 3, 5, 10, 17, 27, 34, 40, 44, 48, 54, 63, 76, 77, 81, 89, 90 \\ 
  \hline
\end{tabular}
\end{table}

The detail solution is attached to this report in an excel file, containing 28 sheets, each sheet for each day.

\subsection{Sensitivity Analysis}

In this section, we will change the parameters to see how variability can affect on our results. 

\subsubsection{Time Horizon}

First, we analyze the impact of the length of time horizon on our results. Table \ref{tab:Timehorizon} shows the results when we increase (decrease) the time horizon. For instance, as we increase the time horizon, it will take more time to find the optimal solution. Besides, our optimal objective value would increase since we have are adding more positive terms to our cost. 

\begin{table}[H]
\caption{Sensitivity Analysis on Time Horizon}
\label{tab:Timehorizon}
    \centering
\begin{tabular}{c p{3.5cm} p{3.5cm} c c} 
 \hline
 T & Id of Opened Buildings & Id of Closed Buildings & Objective Value & CPU Time (S) \\ 
 \hline\hline
14 & 0, 2, 3, 5, 10, 12, 17, 27, 34, 40, 44, 48, 54, 55, 76, 83, 89, 90, (63, 77, 81) & 63, 77, 81, (12, 55, 83) & 21570.1 & 12.43  \\ 
 \hline
 21 & 0, 2, 3, 5, 10, 12, 17, 27, 34, 40, 44, 48, 54, 55, 76, 83, 89, 90, (63, 77, 81) & 63, 77, 81, (12, 55, 83) & 32372.9 & 6.08 \\ 
 \hline
 28 & 0, 2, 3, 5, 10, 12, 17, 27, 34, 40, 44, 48, 54, 55, 76, 83, 89, 90, (63, 77, 81) & 63, 77, 81, (12, 55, 83) & 43175.7 & 12.99 \\ 
 \hline
 35 & 0, 2, 3, 5, 10, 12, 17, 27, 34, 40, 44, 48, 54, 55, 76, 83, 89, 90, (63, 77, 81) & 63, 77, 81, (12, 55, 83) & 53978.5 & 17.02\\ 
 \hline
 42 & 0, 2, 3, 5, 10, 12, 17, 27, 34, 40, 44, 48, 54, 55, 76, 83, 89, 90, (63, 77, 81) & 63, 77, 81, (12, 55, 83) & 64781.3 & 20.40\\ 
 \hline
\end{tabular}
\end{table}

\begin{figure}[H]
 \centering 
 \frame{\includegraphics[width=10cm]{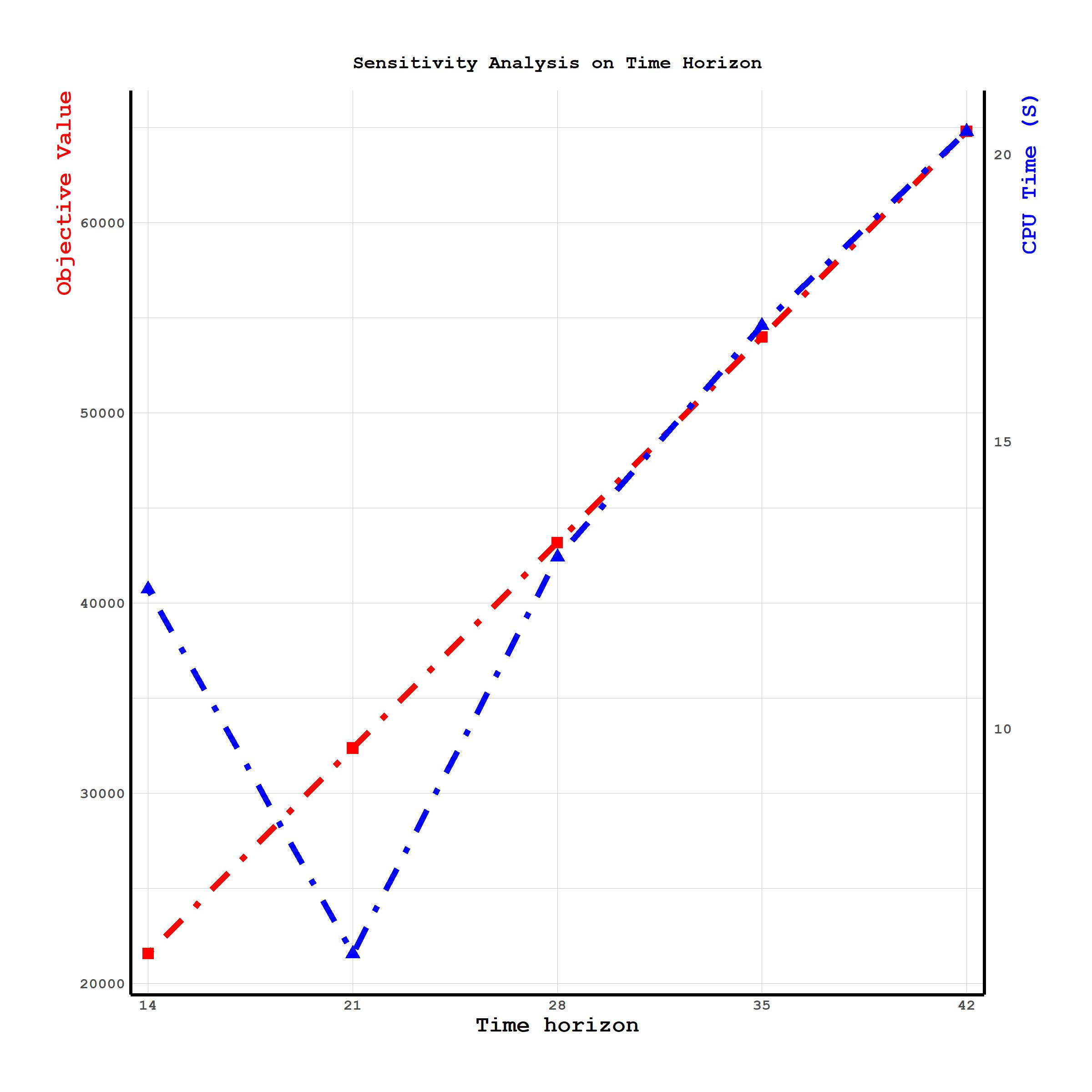}}
 \caption{The impact of the time horizon on the objective functions and CPU process time.}
\label{fig:timehorizon}
\end{figure}

Moreover, as we see the opened and closed buildings remain the same as we increase the time horizon. The reason is that our demand is periodic over the weeks. By increasing the number of weeks, the optimal decision to open and close some stores would be the same. Figure~\ref{fig:timehorizon} summarize the results.

\subsubsection{Number of Facilities to Be Located (P)}

Now we change $P$ and see how it affects our results. Variability of $P$ is important since it would help the company to decide the number of stores they want to buy and invent on the campus.

Note that we could also consider a fixed cost of buying each store in our model (i.e., we could add $Price.P$ to the objective function). However, since in our optimization model $P$ is fixed and given, we don't need to consider it. Figure \ref{fig:P} shows the results when we change $P$. First, it seems that the running time is not directly dependent on $P$. However, $P$ is determining the complexity of the problem. Besides, the objective value is obviously decreasing as we increase $P$ since we didn't consider the price of each store. 

\begin{figure}[H]
 \centering 
 \frame{\includegraphics[width=10cm]{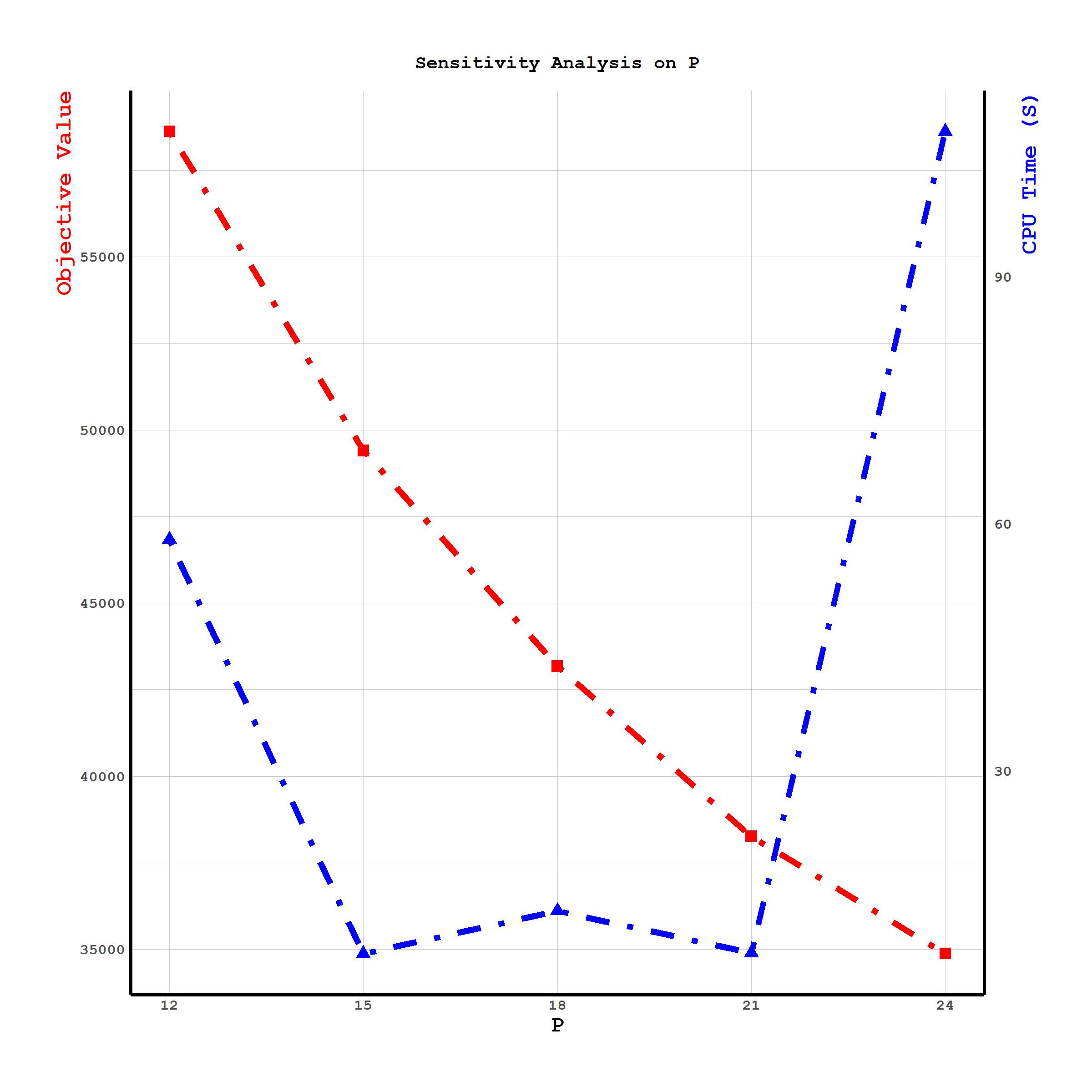}}
 \caption{The impact of the number of facilities on the objective functions and CPU process time.}
\label{fig:timehorizon}
\end{figure}

Also, Table \ref{tab:P} shows how our decision changes in different $P$. There is one interesting point in our decisions. Compare $P=12$ with $P=15$. In $P=15$, we don't use the same opened building we used in $P=12.$ It shows as we want to add one building, we may no longer need another building.

\begin{table}[H]
\caption{Sensitivity Analysis on P}
\label{tab:P}
    \centering
\begin{tabular}{c p{3.5cm} p{3.5cm} c c} 
 \hline
 P & Id of Opened Buildings & Id of Closed Buildings & Objective Value & CPU Time (S) \\ 
 \hline\hline
12 & 2, 5, 10, 17, 27, 35, 44, 54, 65, 76, 83, 90, (63, 89) & 63, 89, (65, 83) & 58616.4  & 58.07 \\ 
 \hline
 15 & 2, 3, 5, 10, 27, 35, 45, 48, 63, 66, 74, 76, 83, 89 , 90 (81, 17, 61, 79) & 81, 17, 61, 79, (83, 3, 66, 76) & 49402.9 & 7.79 \\ 
 \hline
 18 & 0, 2, 3, 5, 10, 12, 17, 27, 34, 40, 44, 48, 54, 55, 76, 83, 89, 90, (63, 77, 81) & 63, 77, 81, (12, 55, 83) & 43175.7 & 12.99 \\ 
 \hline
 21 & 0, 2, 5, 7, 10, 12, 17, 23, 26, 34, 40, 44, 48, 53, 54, 55, 76, 83, 85, 89, 90, (64, 77, 81, 3, 79) & 64, 77, 81, 3, 79, (7, 55, 83, 76, 85) & 38266.5 & 7.89\\ 
 \hline
 24 & 0, 2, 3, 5, 6, 10, 12, 17, 23, 34, 40, 44, 48, 53, 54, 61, 67, 76, 80, 81, 83, 87, 89, 90, (64, 78) & 64, 78, (80, 87) & 34874.4 & 107.60 \\ 
 \hline
\end{tabular}
\end{table}

\subsubsection{Cost Coefficients}
Followed by two previous sections, we want to discuss the effects of changing opening and closing costs together on our results. Table \ref{tab:together} shows the results when we increase (decrease) the cost coefficients. Considering the case that costs are both zero, the model tries to open (and close) any building in each period. In other words, it does not care how many building wants to open or close each day. Another noteworthy point is that running time is increasing as we increase the costs. It shows that the trade-off between keeping current buildings and opening (closing) other buildings is becoming important and challenging to the model. 

\begin{table}[H]
\caption{Sensitivity Analysis on Cost Coefficients (together)}
\label{tab:together}
    \centering
\begin{tabular}{c p{3.5cm} p{3.5cm} c c} 
 \hline
 $\gamma^{o(c)}$ & Id of Opened Buildings & Id of Closed Buildings & Objective Value & CPU Time (S) \\ 
 \hline\hline
 0 & 0, 2, 3, 10, 12, 17, 26, 35, 45, 48, 54, 65, 76, 83, 86, 88, 90, (1, 4, 5, 23, 26, 27, 34, 40, 44, 56, 63, 81, 89) & 4, 5, 27, 34, 40, 44, 56, 63, 81, 89, (23, 26, 35, 45, 65, 86, 88)  & 42900 & 3.90 \\ 
 \hline
 2.5 & 0, 2, 3, 5, 10, 12, 17, 27, 34, 40, 44, 48, 54, 55, 76, 83, 89, 90, (63, 77, 79, 81) & 63, 77, 79, 81, (12, 55, 76, 83) & 43049.3 & 9.34 \\ 
 \hline
 5 & 0, 2, 3, 5, 10, 12, 17, 27, 34, 40, 44, 48, 54, 55, 76, 83, 89, 90, (63, 77, 81) & 63, 77, 81, (12, 55, 83) & 43175.7 & 12.99  \\  
 \hline
 10 &  0, 2, 3, 5, 10, 12, 17, 27, 34, 40, 44, 48, 54, 55, 76, 83, 89, 90, (80, 81) & 80, 81, (12, 83)  & 43382.3 & 58.45 \\ 
 \hline
 20 &  0, 2, 3, 5, 10, 12, 17, 27, 34, 40, 44, 48, 54, 55, 76, 83, 89, 90, (81) & 81, (83)  & 43658.2 & 82.58 \\ 
 \hline
\end{tabular}
\end{table}

\begin{figure}[H]
 \centering 
 \frame{\includegraphics[width=10cm]{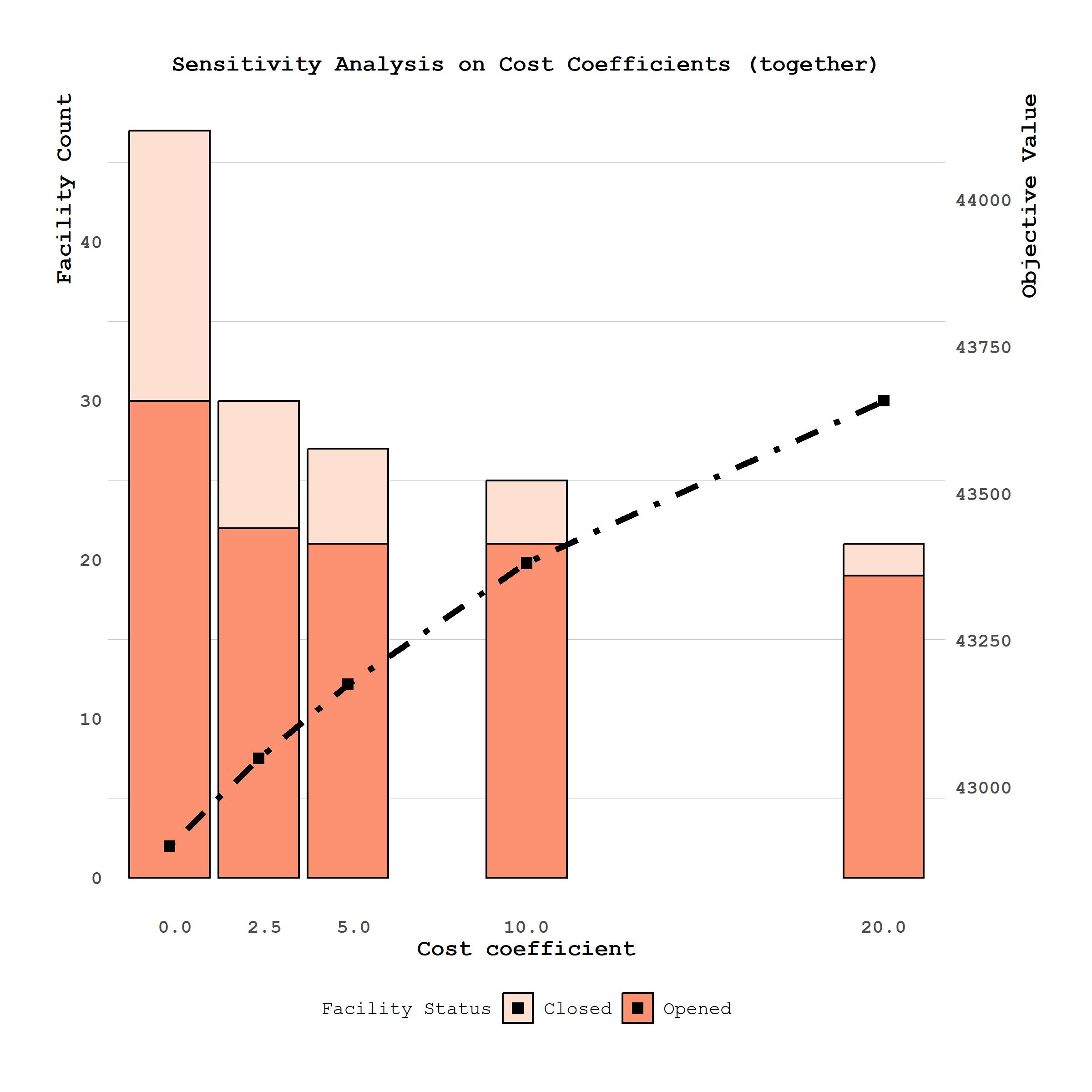}}
 \caption{The impact of opening/closing cost of facilities on the objective functions.}
\label{fig:together}
\end{figure}

Figure~\ref{fig:together} shows that the initially opened buildings remain the same and starts to dynamically be changed as we increase the opening/closing cost. Moreover, the objective value increase as we increase the opening cost. Also, for the large opening costs, we prefer not to open new buildings (and close the current open buildings) since it would be costly. 

\section{Discussion}\label{sec:conclusion}
This study has presented a mixed-integer linear formulation for the Dynamic Modified Stochastic p-Median Problem in a Competitive Supply Chain Network Design. The proposed model takes into account the robust optimization and time horizon as a novel approach, which enables the decision-maker to consider uncertainty and short-term changes in the supply chain network design. The robust optimization approach used in this study allows for the consideration of different scenarios and uncertainty in demand and supply, which is crucial in real-world applications. Additionally, the time horizon approach allows for the dynamic nature of the problem to be captured, which is important in today's fast-paced business environment.

The proposed model was tested using computational experiments, and the results demonstrate the effectiveness of the proposed approach in handling the dynamic and stochastic nature of the problem. The results also provide valuable insights for practitioners and researchers in the field of supply chain network design. The proposed model can be extended and applied to other similar problems in the field, such as facility location, transportation and logistics, and inventory management.

One of the main contributions of this study is the integration of robust optimization and time horizon in the mixed-integer linear formulation for the Dynamic Modified Stochastic p-Median Problem. The robust optimization approach allows for the consideration of different scenarios and uncertainty, while the time horizon approach allows for the dynamic nature of the problem to be captured. This integration provides a more realistic and practical solution to the problem, which can be useful for practitioners and researchers in the field.

Another important contribution of this study is the application of the proposed model to a competitive supply chain network design problem. This application is relevant and valuable as it provides insights into how the proposed model can be used in a real-world context. The results of the computational experiments demonstrate the effectiveness of the proposed model in handling the dynamic and stochastic nature of the problem and provide valuable insights for practitioners and researchers in the field. For instance, we can change over the problem definition to solve any other location problems. Covering location problems are valuable to be investigated, such as trying to find the optimal number of Robomarts \cite{pymnts.com_2018} that can serve all facilities on the workspace if the single Robomart can only serve facilities within limited miles.

This study has the potential usefulness of the proposed model in the context of a pandemic and quarantine. The ability to handle uncertainty and short-term changes in the supply chain network design. The robust optimization approach used in this study allows for the consideration of different scenarios and uncertainty in the demand and supply, which is crucial in a pandemic context. The proposed model can be used to design and optimize supply chain networks that are more resilient and adaptable to the changing conditions caused by a pandemic. This can help businesses and organizations to minimize disruptions and maintain the continuity of operations during challenging times.

One limitation of this study is that the Robomart can only be located at specific facilities (nodes) in the proposed model. In reality, however, the Robomart can also be located somewhere on the route between two nodes (edges). This limitation may affect the validity and applicability of the proposed model in real-world scenarios. To address this limitation, it is possible to consider analogous absolute $p$-median problems to simulate reality. This approach would involve the inclusion of edge-based locations for the Robomart in the model, which would provide a more realistic representation of the problem. However, this would require additional mathematical development and computational resources, and would be a subject for future research.

Another potential area of future research is to make the problem an Adaptive Robust Optimization (ARO) problem \cite{sun2021robust}. In this approach, the number of trucks $p$ would be determined in the first stage and other variables in the second stage. This would enable a more flexible and dynamic approach to supply chain network design, as the number of trucks can be adjusted in response to changes in demand and supply.

In summary, the proposed DMS-p-MP model takes into account the robust optimization and time horizon as a novel approach, which enables the decision-maker to consider uncertainty and short-term changes in the supply chain network design. The results of the computational experiments demonstrate the effectiveness of the proposed model in handling the dynamic and stochastic nature of the problem and provide valuable insights for practitioners and researchers in the field. The proposed model can be extended and applied to other similar problems in the field of supply chain network design. The integration of robust optimization and time horizon in the proposed model provides a more realistic and practical solution to the problem, which can be useful for practitioners and researchers in the field.

\bibliographystyle{elsarticle-harv}\biboptions{authoryear}
\bibliography{cas-refs}

\end{document}